\documentclass[12pt]{amsart}

   \usepackage{latexsym}
   \usepackage{amsmath}
   \usepackage{amssymb}
   \usepackage{amscd}
   \usepackage{multicol}
\usepackage{xypic}	

   \newtheorem{theorem}{\bf Theorem}[section]
   \newtheorem{theo}{\bf Theorem}[subsection]
   \newtheorem{corollary}[theorem]{\bf Corollary}
   
   \newtheorem{lemma}[theorem]{\bf Lemma}
   \newtheorem{definition}[theorem]{\bf Definition}

   \newcommand{\res}{ {\rm Res} }

\newcommand{\ih}{I\!H}
\newcommand{\quott}{/\!/}

   \def\R{\mathbb{R}}
   \def\C{\mathbb{C}}
   \def\Z{\mathbb{Z}}

   \def\t{\mathfrak{t}}
   
    \def\k{\mathfrak{k}}
    
   \def\quott{/\! /}
   \def\norm{|\! |}
   
   \def\t{\mathfrak{t}}

   \def\k{\mathfrak{k}}
   
   \def\s{\mathfrak{s}}
   
   \def\cF{\mathcal{F}}
   \def\eps{\epsilon}
   \def\D{\mathcal{D}}

 \renewcommand{\ker}{ {\rm Ker}}

\include{include}

\begin{document}

\noindent {\LARGE \bf The Kirwan map, equivariant Kirwan maps, and their
 kernels}
\bigskip\\
{\bf Lisa C. Jeffrey } \\
Department of Mathematics, University of Toronto,
Toronto ON, M5S 3G3, Canada \\
E-mail address: {\tt jeffrey@math.toronto.edu }
\smallskip \\
{\bf  Augustin-Liviu Mare}\\
{Department of Mathematics and Statistics, University of Regina,}
   {Regina SK,  S4S 0A2, Canada}, 
    E-mail address: {\tt mareal@math.uregina.ca}
\smallskip \\
{\bf Jonathan M. Woolf}\\
{Christ's College, Cambridge, CB2 3BU, UK}, 
   E-mail address: {\tt jw301@cam.ac.uk}
\bigskip
\title{}

\vspace{-1.5cm}

\begin{abstract}
Consider a Hamiltonian  action
    of a compact  Lie group $K$ on a compact symplectic
   manifold. We find  descriptions of the kernel of the Kirwan map
 corresponding to a regular
   value of the moment map $\kappa_K$. We start with  the case   when  $K$
is
 a torus $T$:
   we determine the
   kernel of the equivariant Kirwan map
   (defined by Goldin in [Go]) corresponding to
   a generic circle $S\subset T$, and  show
   how to recover from this the kernel of $\kappa_T$,
   as described by  Tolman and Weitsman in [To-We].
   (In the situation when the fixed point set of the
   torus action is finite,
   similar results have been obtained in our previous papers [Je],
[Je-Ma]).
   For a compact nonabelian
   Lie group $K$ we will use the ``non-abelian localization formula'' of
   [Je-Ki1] and [Je-Ki2]  to establish  relationships --- some of them
 obtained by
   Tolman and Weitsman in [To-We] ---  between $\ker(\kappa_K)$ and
 $\ker(\kappa_T)$,
   where $T\subset K$ is a maximal torus. In the appendix we prove that the same 
   relationships remain true in the case when $0$ is no longer a regular value
of $\mu_T$.

\end{abstract}

\maketitle




\section{Introduction}

   Let $M$ be a compact symplectic manifold
   acted\footnote{The fixed point set of the action
    is not assumed to consist of isolated points, unlike the hypothesis in
 [Je] and [Je-Ma]).} on by a torus $T$
   in a Hamiltonian fashion.  Assume that $0$ is a regular
   value of the moment map $\mu:M\to \t^*=({\rm Lie} \ T)^*$ and consider
   the {\it symplectic reduction}
   $$M\quott T=\mu^{-1}(0)/T.$$
   By a   well-known result (see \cite{Ki}), the Kirwan map $\kappa
 :H_T^*(M)\to
   H^*(M\quott T)$ induced by the inclusion
   $\mu^{-1}(0)\hookrightarrow M$ via the identification $H_T^*(M)=
   H^*(M\quott T)$ is surjective\footnote{In this paper all cohomology
groups
 are with complex coefficients.}.
   Determining  the kernel of
   $\kappa$ is an important problem: combined with a presentation of
   the equivariant cohomology ring $H_T^*(M)$, this gives a
   description of the cohomology ring of $M\quott T$.

   Now consider a 1-dimensional torus $S\subset T$ with Lie algebra
   $\s$ and let $\mu_S:M \to \s^*$ be the  moment map corresponding
   to the action of $S$ on $M$. We will call $S$ {\it generic} if

   (i) the two fixed point sets $M^S$ and $M^T$ are equal

   (ii) $0$ is a regular value of the moment map $\mu_S:M\to \s^*$.

   Denote by $M\quott S=\mu_S^{-1}(0)/S$ the symplectic reduction
   corresponding to the $S$ action on $M$. In the same way as before, we
   can consider the map induced by the inclusion
 $\mu_S^{-1}(0)\hookrightarrow M$
   at the level of the $T$-equivariant cohomology and identify at the
   same time
   $$H_T^*(\mu_S^{-1}(0))=H_{T/S}^*(\mu_S^{-1}(0)/S).$$
   Goldin \cite{Go} called the resulting map
   $$\kappa_S:H_T^*(M)\to H_{T/S}^*(M\quott S)$$
   the {\it equivariant Kirwan map} and she proved:

   \begin{theorem}\label{rebecca}{\rm (see \cite{Go})} If $0$ is  a regular
 value of
   $\mu_S$, then the  map $\kappa_S$ is surjective.
   \end{theorem}

   Our first theorem (Theorem \ref{firstmain})
 shows that the knowledge of $\ker(\kappa_S)$, for a
 generic circle $S\subset T$,
    is crucial for determining the kernel of the Kirwan map
   $$\kappa:H^*_T(M)\to H^*(M\quott T).$$
   \begin{theorem}\label{firstmain}
   $$\ker (\kappa) = \sum_{S} \ker (\kappa_{S}),$$
   where the sum runs over all generic circles $S\subset T$.
   \end{theorem}

   One aim  of our paper is to
    describe the kernel of $\kappa_{S}$. Let us denote
   first by $\cF$ the set of all connected components of the fixed
   point set $M^T$: note that the moment map $\mu$ is constant on
   each $F\in \cF$. To any generic circle $S\subset T$ we assign the
   partition of $\mathcal{ F}$ given by ${\mathcal F}={\mathcal
   F}_-\cup {\mathcal F}_+$, where
   $${\mathcal F}_+=\{F\in {\mathcal F}: \mu_S(F)> 0\}$$
   and
   $${\mathcal F}_-=\{F\in {\mathcal F}: \mu_S(F)<0\}.$$
    Also, let $X,Y_1,\ldots,Y_m$ be variables corresponding to
   an integral basis of $\t^*$ such that
   \begin{equation}\label{basis}X|_{\s} \ {\rm corresponds \ to \ an \
   integral \ element \ of} \ \s^* {\rm \ and} \ Y_j|_{\s}=0, j\ge
   1.\end{equation} In this way, we have
   $$H_T^*({\rm pt})=S(\t^*)=\C[X,Y_1,\ldots,Y_m].$$
  The following definition (Definition
\ref{reskerdef} for ${\rm Ker}_{\rm res}$)\footnote{For more details, we refer the reader to
 section 2.}
is suggested by the residue formula of
    [Je-Ki1] and
   [Je-Ki2]:
   \begin{equation}\label{residue}\ker(\kappa_S)=\ker_{{\rm
 res}}(\kappa_S).\end{equation}

   \begin{definition} \label{reskerdef}
    \noindent (i)
   The {\rm residue
   kernel} $\ker_{{\rm res}}(\kappa_S)$ is the set of all classes $\eta\in
 H_T^*(M)$ with the
   property that
   $${\rm Res}_X^+\sum_{F\in\cF_+}\int_F\frac{i_F^*(\eta\zeta)}{e_F} =0$$
   for all $\zeta\in H_T^*(M)$. Here $e_F=e(\nu F)\in H_T^*(F)$ is the
   equivariant Euler class of the normal bundle of $F$ and
   ${\rm Res}_X^+$ is defined by
   \begin{equation}\label{sum}{\rm Res}_X^+(h)=\sum_{b\in\C}{\rm
 Res}_{X=b}h\end{equation}
   for every  rational function $h$ in the variables $X,Y_1,\ldots,Y_m$.

   \noindent (ii) We define the residue kernel of $\kappa $ as
   $$\ker_{\rm res} (\kappa) = \sum_S \ker_{\rm res} (\kappa_S) $$
    \end{definition}

   We shall prove the following result:
    \begin{theorem}\label{secondmain} If $S\subset T$ is a generic circle,
 then we have that
\begin{equation} \label{e:secondmain}
\ker_{{\rm res}} (\kappa_S) = K_-^S \oplus K_+^S, \end{equation}
   where $K_-^S$ (resp. $K_+^S$) denote the set of all equivariant cohomology
   classes $\eta$ whose restriction to ${\mathcal F}_-^S$ (resp.
   ${\mathcal F}_+^S$) is zero.\end{theorem}

We remark  that  we prove Theorem
\ref{secondmain} directly, without using the fact that the
right and left hand sides of (\ref{e:secondmain})
are known to be the kernel of $\kappa_S$ (by results
of \cite{Je-Ki1} and Goldin \cite{Go}).

   As a corollary of
Theorem \ref{secondmain}, we deduce the following description
of
 $\ker_{\rm res}(\kappa)$:
   \begin{corollary}{\label{corollary}} The
residue kernel of
the
 Kirwan map
   $\kappa:H_T^*(M)\to H^*(M\quott T)$ is given by
   $$\ker_{\rm res} (\kappa) = \sum_S (K_-^S \oplus K_+^S),$$
   where $S$ is as in Theorem \ref{firstmain}.
   Here $K_{\pm}^S$ consist of all equivariant cohomology
   classes $\alpha$ which restrict to zero on all components
   $F\in {\mathcal{F}}$ with the property that
    $$\pm\mu(F)(\xi)  >  0 $$
   where  $\xi$ is a fixed non-zero vector in the Lie algebra
   of  $S$.
   \end{corollary}
 \noindent  {\em Proof of Corollary:}
  The corollary follows immediately from Theorem \ref{secondmain}
   because Theorem \ref{secondmain} tells us that
   $${\rm Ker}_{\rm res} (\kappa_S) = K_-^S \oplus K_+^S$$
   so
   $$\sum_S {\rm Ker}_{\rm res} (\kappa_S) = \sum_S ( K_-^S \oplus K_+^S)$$
   which is ${\rm Ker}_{\rm res}(\kappa)$ by definition
 (Definition \ref{reskerdef}). \hfill $ \square$

Tolman-Weitsman's theorem \cite{To-We} is as follows:
\begin{theorem} \label{t:tw} We have
 \begin{equation} \label{e:tw}
 {\rm Ker} (\kappa) = \sum_S ( K_-^S \oplus K_+^S).  \end{equation}
\end{theorem}

 \noindent{\em Proof of Tolman-Weitsman's theorem:}
We obtain Tolman-Weitsman's theorem
by the following steps.

\begin{description}
\item[Step 1] Theorem \ref{firstmain} tells us that
${\rm Ker}(\kappa) = \sum_S {\rm Ker}(\kappa_S)$
\item[Step 2] Results of \cite{Je-Ki1} and \cite{Je-Ki2} tell us
that  ${\rm Ker} (\kappa_S) = {\rm Ker}_{\rm res} (\kappa_S)$
\item[Step 3] By Theorem \ref{secondmain}
${\rm Ker}_{\rm res} (\kappa_S) = K_-^S \oplus K_+^S. $ \hfill $\square$
\end{description}
The key step is Step 3, which uses residues which enable us to easily
exhibit the
structure of the residue kernel as a sum of classes vanishing on one side of
certain hyperplanes.
   We must emphasize that it is unexpected that one can
 deduce Tolman-Weitsman's result (Theorem \ref{t:tw})
 by this means; obtaining
a new
   proof of this result was not one of our original goals, but rather an
 unanticipated byproduct of our analysis.

 Our goal was to prove directly that the residue kernel
${\rm Ker}_{\rm res}$ (see Definition \ref{reskerdef})
 equals the Tolman-Weitsman kernel
 (i.e. the right hand side of (\ref{e:tw})).
 We obtain this result (Corollary \ref{corollary}) by combining
 Theorem \ref{firstmain} and Theorem \ref{secondmain}.
 This directly generalizes the main result of
 \cite{Je}. The two descriptions of the kernel of the
 Kirwan map are apparently quite different, so it is
 illuminating to see directly that they coincide.

   The final goal of the paper is to describe the kernel of the
   Kirwan map in the case of the Hamiltonian action of a non-abelian
   Lie group. Let $K$ be an arbitrary compact connected Lie group with
 maximal
   torus $T\subset K$ and let $W=N_K(T)/T$ be the corresponding Weyl
   group. If $M$ is a $K$ manifold, then there exists a natural action
   of $W$ on $H_T^*(M)$. If $H^*_T(M)^W$ denotes the space of $W$-invariant
   cohomology classes, then we have the isomorphism
   $$H^*_K(M)=H^*_T(M)^W,$$
   (see [At-Bo1]) which allows us to identify the two spaces. Let us also
 denote by $\D$ the element of $H^*_T({\rm pt})=S(\t^*)$
   given by the product of all positive roots.  Suppose now that $M$ is a
 compact symplectic
   manifold and the action of $K$ is Hamiltonian. Assume that
   $0$ is a regular value of the moment map
   $$\mu_K:M\to \k^*$$
   and consider
   the symplectic quotient
   $$M\quott K =\mu_K^{-1}(0)/K.$$
   The Kirwan surjection (see [Ki])
   $$\kappa_K:H^*_K(M)\to H^*(M\quott K)$$
   is, exactly as in the abelian case, induced by restriction from $M$ to
 $\mu_K^{-1}(0)$ followed by the isomorphism
   $$H^*_K(\mu_K^{-1}(0))\simeq H^*(\mu_K^{-1}(0)/K).$$

   The following result was proved in [Je-Ki1]:

   \begin{theorem}{\rm (non-abelian residue formula [Je-Ki1], [Je-Ki2]).}
   For every  cohomology class $\eta\in H^*_K(M)=H^*_T(M)^W$ we have that
   \begin{equation}{\label{nonabelianresidue}}\kappa_K(\eta)[M\quott
 K]=c_1{\rm Res}(\sum_{F\in \cF}\int_F\frac{\D^2\eta|_F}{e_F})
   \end{equation}
   where $c_1$ is a non-zero constant, $\cF$ is the set of connected
 components of the fixed point set $M^T$  and
   {\rm Res} is a certain complex valued operator\footnote{For
   the exact ``residue formula'', the reader is referred to [Je-Ki1] and
 [Je-Ki2];
   but for the purposes of   the present paper, equation
 (\ref{nonabelianresidue}) is
   sufficient. A more detailed description is
 given in  Section 2.2, where it is used.}
   on the space of  functions  of the form $p/q$, with $p,q\in S(\t^*)$,
 $q\ne 0$.
   \end{theorem}
   In the case $K=T$, the class $\D$ is $1$, and we obtain
   immediately the {\it abelian residue formula}: for every $\eta\in
H^*_T(M)$
 we have
    that
   \begin{equation}{\label{abelianresidue}}\kappa_T(\eta)[M\quott T]=
   c_2{\rm Res}(\sum_{F\in \cF}\int_F\frac{\eta|_F}{e_F})
   \end{equation}
   where ${\rm Res}$ is the same operator as the one from equation
 (\ref{nonabelianresidue})
   and $c_2$ is a non-zero constant.
   Also note that the residue formula (\ref{nonabelianresidue}) gives a
 complete  description
   of $\ker(\kappa_K)$, since for $\eta\in H^*_K(M)$ we have that
   \begin{equation}\label{equivalence}\kappa_K(\eta)=0 \ \quad \ {\rm iff }
 \quad \kappa_K(\eta\zeta)[M\quott K]=0\ {\rm for\  every \ }
   \zeta\in H^*_K(M).\end{equation} We would like to make this description
 more explicit, by giving
   a direct relationship between
   $\ker(\kappa_K)$ and $\ker(\kappa_T)$ (the latter being given by Corollary
 \ref{corollary}).
   We will prove that:
   \begin{theorem}\label{nonabelian} Take $\eta\in H^*_K(M)= H^*_T(M)^W$.
   The following assertions are equivalent:
   \begin{itemize}
   \item[(i)] $\kappa_K(\eta)=0$,
   \item[(ii)] $\kappa_T(\D\eta)=0$,
   \item[(iii)] $\kappa_T(\D^2\eta)=0$.
   \end{itemize}
   \end{theorem}

\noindent   The proof will be given in section 5. Some explicit relationships
between
 $\ker(\kappa_K)$ and $\ker(\kappa_T)$ will also
   be discussed there.

 The layout of our paper is as follows. In
 Section 2 we recall the characterizations of residues  from
 \cite{Gu-K},
 \cite{Je-Ki1}, \cite{Je-Ki2}, and
summarize
 some results from Morse theory which will be used later.
  Section 3 is devoted to a proof of Theorem \ref{secondmain},
 and Section 4 to the proof of Theorem \ref{firstmain}.
 Finally Section 5 gives a characterization of the kernel
 of $\kappa$ for nonabelian group actions.

   {\bf Remarks.}
 \begin{description}
\item[1]  The chief advances described in this article  are as follows.
\begin{description}
\item[(a)] We prove Theorem \ref{firstmain} (which relates the kernel of the
Kirwan map to the kernels of the
equivariant Kirwan maps $\kappa_S$).
\item[(b)] We prove Theorem \ref{secondmain} (which characterizes the residue
kernel ${\rm Ker}_{\rm res} (\kappa_S)$ as
linear combinations of elements of equivariant cohomology
vanishing on the preimage of one side of a hyperplane).
\item[(c)] Our construction does not assume isolated fixed points.
 Theorem \ref{secondmain} is a generalization of the main results of
 \cite{Je} and \cite{Je-Ma},
   where  actions with isolated fixed points were considered.
\item[(d)] We give an independent characterization
(Theorem \ref{nonabelian}) of the kernel of the Kirwan
map for nonabelian group actions.
 According to Tolman and Weitsman [To-We], the equivalence
   (i) $\iff$(ii) from Theorem
   \ref{nonabelian} can be deduced from S. Martin's integral formula (see [Ma, Theorem B]). 
Our  proof of this equivalence
relies on   the  residue formulas (\ref{nonabelianresidue}) and
   (\ref{abelianresidue}).
\end{description}
   \item[2]  The characterization of $\ker(\kappa_S)$  given in Theorem
 \ref{secondmain} (see  equation (\ref{residue})) was proved by
   Goldin in [Go], by using the theorem of Tolman and Weitsman [Theorem 3,
 To-We]. The strategy of our paper
   is different: we first prove  the formula  stated in Theorem
 \ref{secondmain} and then deduce the theorem of Tolman and
   Weitsman starting from Corollary \ref{corollary}  using Theorem
\ref{firstmain}
and the known fact that ${\rm Ker} (\kappa_S) = {\rm Ker}_{\rm res}
(\kappa_S). $

\end{description}

   \section{Residue formulas and  Morse-Kirwan theory}

   The goal of this section is to provide some fundamental
   definitions and results --- which will be needed later ---
   concerning the notions mentioned in the title.

   \subsection{The Kirwan map in terms of residues}
   Let us consider the compact symplectic manifold $M$ equipped with
   the Hamiltonian action of the torus $T$. Let $S\subset T$ be a
   generic circle and consider   the variables $X,Y_1,\ldots,Y_m$
   determined by (\ref{basis}) (see section 1). The following result
   was proved in [Je-Ki1] and [Je-Ki2]:

   \begin{theo}\label{jk} For all $\eta\in H_T^*(M)$ we
   have\footnote{Note that both sides of the equation are in $\C[Y_1,
   \ldots, Y_m]$.}
   \begin{equation}\label{jkformula} \kappa_S(\eta)[M\quott S]= c\sum_{F\in
 \cF_+}{\rm Res}_X^+
   \int_F\frac{i_F^*(\eta)}{e_F}.\end{equation}  Here $e_F\in H_T^*(F)$
 denotes
   the $T$-equivariant Euler class of the normal bundle $\nu(F)$, $c$
   is a non-zero constant\footnote{The precise value of $c$ is not
   needed in our paper, but the interested reader can find it in
   [Je-Ki2, \S 3].}, and the meaning of {${\rm Res}_X^+\int_F$} is
   given in (\ref{e:resplus})
   \end{theo}

   For all $\alpha\in H_T^*(F)$, we obtain
   $${\rm Res}_X^+\int_F\frac{\alpha}{e_F}$$
   as follows: We may assume that the normal bundle $\nu(F)$  splits
   as$$\nu(F)=L_1\oplus \ldots \oplus L_k$$ where $L_i$, $1\le i \le
   k$ are $T$-equivariant complex line bundles. We express the
   reciprocal of $e_F$  as
   \begin{equation}\label{reciprocal}\frac{1}{e_F}=
   \prod_{i=1}^{k}\frac{1}{(m_iX+\beta_i(Y)+c_1(L_i))}=
   \prod_{i=1}^k\sum_{r_j\ge
   0}\frac{(-c_1(L_i))^{r_j}}{(m_iX+\beta_i(Y))^{r_j+1}}.\end{equation}
   Here $m_i$ is the weight of the representation of $S$ on $L_i$,
   $\beta_i(Y)$ is a certain linear combination of $Y_1,\ldots, Y_m$  with
 integer
   coefficients,  and $c_1(L_i)$ is the first
   Chern class of $L_i$. Since $c_1(L_i)$ is in $H^2(F)$, it is
   nilpotent, hence all sums involved in (\ref{reciprocal}) are
   finite. The class $\alpha$ is in $H_T^*(F)=H^*(F)\otimes
   \C[X,\{Y_i\}]$, hence by (\ref{reciprocal}),
   $$\int_F\frac{\alpha}{e_F}$$ appears as a sum of
    rational functions of the type
   $$h(X,Y_1,\ldots,Y_m)=\frac{p(X,Y_1,\ldots, Y_m)}{q(X,Y_1,\ldots,
Y_m)}$$
   where $p,q\in \C[X,\{Y_i\}]$, with $q=\prod_j(n_jX+\sum_l
   n_{jl}Y_l) $,  $n_j\in \Z\setminus\{0\}$ and $n_{jl}\in\Z$. In
   order to define ${\rm Res}_X^+(h)$, regard $Y_1,\ldots,Y_m $ as
   complex constants  and set
    \begin{equation} \label{e:resplus}
 \text{Res}_X^+(h)=\sum_{b\in\C}\text{Res}_{X=b}\frac{p}{q} \end{equation}
   where on the right hand side $\frac{p}{q}$ is interpreted as a
   meromorphic function in the variable $X$ on $\C$.

   {\bf Remark 2.1.2} Guillemin and Kalkman gave another definition
   of the residue involved in Theorem \ref{jk}.
    (see  [Gu-Ka, section 3]).
    Write

\begin{align*}\frac{1}{e_F}&=\prod_{i=1}^k\frac{1}{(m_iX+\beta_i(Y)+c_1(L_i)
 )}\\
 {}&=\prod_{i=1}^k\frac{1}{m_iX(1+\frac{\beta_i(Y)+c_1(L_i)}{m_iX})}\\{}
 &=\prod_{i=1}^k
(\frac{1}{m_iX}-\frac{\beta_i(Y)+c_1(L_i)}{(m_iX)^2}+\frac{(\beta_i(Y)+c_1(L
 _i))^2}{(m_iX)^3}
   -\ldots ).\end{align*}
If we multiply the $k$ power series in $X$
in the last expression and then multiply the result by the
polynomial $\alpha\in H^*(F)\otimes\C[X,\{Y_i\}]$, we obtain a
series
$$\sum_{r=r_0}^{-\infty}\gamma_rX^r,$$
where $\gamma_r$ are in  $H^*(F)\otimes \C[\{Y_i\}]$. The
Guillemin-Kalkman residue is $\int_F\gamma_{-1}$.  It is a simple
exercise to show that the latter coincides with
the expression $${\res}_X^+\int_F\frac{\alpha}{e_F}$$ of [Je-Ki2]
which we have
defined above.

 \newcommand{\Jac}{\bigtriangleup}

   \begin{definition} \label{e:resdef} (see Section 3 of \cite{Je-Ki2})
 If $h$ is a meromorphic function of $m+1$ variables $X_1, \dots, X_{m+1}$
 then
   \begin{equation}
   \res(h(X))=  \Jac \res^+_{X_1} \circ  \ldots \circ
   \res^+_{X_{m+1}}  (  h(X) )\end{equation}
   where the variables $X_1,\ldots,X_{m}$ are held constant while
   calculating $\res^+_{X_{m+1}}$, and $\Jac$ is the determinant of some
   $(m+1)\times (m+1)$ matrix whose columns are the coordinates of an
   orthonormal basis of ${\bf  t} $ defining the same orientation as the
   chosen coordinate system.
 Here  the symbols $\res^+_{X_j} $ were defined in (\ref{e:resplus}).
   \end{definition}

 \begin{definition}
 If $\alpha, \beta \in H^*_T(M)$ then
 we define the complex number
 $${\rm Res}(\alpha \beta) =
 {\rm Res}(\sum_F \int_F \frac{\alpha \beta}{e_F}). $$
 Here $F $ are the components of the fixed point set of the $T$
 action, and $e_F$ is the equivariant Euler class of the normal
 bundle.
 \end{definition}

By (\ref{residue})  and
Theorem \ref{firstmain} (which will be proved in
Section 4)   we have
\begin{equation} \label{e:kerres}
{\rm Ker} (\kappa) = {\rm Ker_{res} } (\kappa).
\end{equation}
By Corollary \ref{corollary}  we have also  that
$$ {\rm Ker} (\kappa) = \sum_S (K_-^S \oplus K_+^S). $$
By \cite{Je-Ki1} we have that
\begin{equation}{\rm Ker}(\kappa) = \{ \alpha \in H^*_T(M) |
{\rm Res} (\alpha \beta) = 0 ~\forall ~ \beta \in
H^*_T(M) \}.  \end{equation}
Combining these observations we obtain
\begin{theorem} \label{t:restw}
A class $ \alpha \in H^*_T(M) $
satisfies ${\rm Res}(\alpha \beta) = 0 $ for all $\beta \in H^*_T(M)$ if and
only if $ \alpha \in \sum_S K^-_S \oplus K^+_S$.
\end{theorem}
\hfill  $\square$

   \subsection{Kirwan-Morse theory}
   Let $S\subset T$ be  a generic circular subgroup of a torus
   $T$ which acts in a Hamiltonian manner on the compact symplectic
   manifold $M$.  One knows that the critical points of the moment
   map
   $$\mu_S=f:M\to \s^*=\R$$
   are the fixed points of the $S$ action: as usual, we will denote
   by $\cF$ the set of connected components of $$M^T=M^S={\rm
   Crit}(f).$$ Moreover, for  every $F\in \cF$, the restriction of the
   Hessian of $f$ to the normal space  to $F$ is nondegenerate, which
   means that $f$ is a Morse-Bott function. We assume for simplicity
   that the values of $f$ on any two elements of $\cF$ are different.

   Fix $F\in \cF$ and let $\nu F$ be the normal bundle.  Consider the
   $T$-invariant splitting, $$\nu F=\nu^- F \oplus \nu^+ F$$ into the
   negative and positive subspaces of the Hessian.  The rank of the
   bundle $\nu^-F$ is just the index ${\rm Ind}(F)$ of the critical
   level $F$. The following result follows from [At-Bo2] (for more
   details, see [Go, section 2]):
   \begin{theo}\label{morse} Suppose that $\eps$ is sufficiently small that
   $f(F)$ is the only critical value in the interval
   $(f(F)-\eps,f(F)+\eps)$ and consider
   $$M_-=f^{-1}(-\infty, f(F)-\eps), \quad M_+=f^{-1}(-\infty,
   f(F)+\eps).$$

   a) There exists a $H_T^*({\rm pt})$-linear map\footnote{The map $r$ is
in
 fact the restriction map
   $H^*_T(M_+,M_-)\to H^*_T(M_+)$ composed with the inverse of the excision
 isomorphism
   $H^*_T(M_+,M_-) \to H^*_T(D^{{\rm Ind}F}, S^{{\rm Ind}F-1})$ and then
with
 the inverse
   of the Thom isomorphism $H^*_T(D^{{\rm Ind}F}, S^{{\rm Ind}F-1})\to
 H^{*-{\rm Ind}F}_T(F)$,
   where $D^{{\rm Ind}F}$ is  the  cell bundle over $F$ attached to $M_-$
 after passing the critical level
   $F$, and  $S^{{\rm Ind}F-1}$ is the corresponding sphere bundle over
 $F$. }
   $$r:H_T^{*-{{\rm Ind} F}}(F)\to H_T^*(M_+)$$
   such that
   
   $\ $ $\bullet$ the sequence
   \begin{equation}\label{exactsequence}0\to H_T^{*-{{\rm Ind}
 F}}(F)\stackrel{r}{\to} H_T^*(M_+)\to H_T^*(M_-)\to 0\end{equation}
   is exact,

   $\ $ $\bullet$ the composition of  the restriction map
   $$ H_T^*(M_+)\to
   H_T^*(F)$$ with $r$ is the multiplication by the equivariant Euler
   class $e_T(\nu^-(F))$; moreover, the image of $r$ consists of all
classes
 in
   $H^*_T(M_+)$ whose restriction to $F$ is a multiple of $e_T(\nu_f^-F)$.

   b) The maps $H_T^{*}(F)\to H_T^*(F)$ given by the multiplication
   by $e_T(\nu^-(F))$ and $e_T(\nu^+(F))$ are injective (in other
   words $e_T(\nu^{\pm}(F))$ are not  zero divisors).

   c) The restriction map $$H_T^*(M)\to H_T^*(M_+)$$ is surjective.

   d)  There exists a class
   $\alpha^+(F)\in H_T^*(M)$ with the following properties:

   $ \  \bullet$ $\alpha^+(F)|_G=0$ for all $G\in \cF$
   with $f(G)> f(F)$).

   $ \  \bullet$ $\alpha^+(F)|_F=e_T(\nu^+(F))$

   \end{theo}

   {\bf Remark 2.2.2} A similar  result holds if we take instead of
$f$ (a generic component of the moment map $\mu:M\to \t^*$) the
function
$$g(x)=\norm \mu(x)\norm^2, \quad x\in M.$$
 Although $g$ is no longer Morse-Bott, like $f$, it is {\it
minimally degenerate} in the sense defined by Kirwan [Ki]. In
particular the critical set Crit($g$) is a disjoint union of
$T$-invariant closed subspaces $C$, call them the {\it critical
sets of $g$}, with the following properties:
\begin{itemize}
\item $\mu$ is constant on $C$;
\item the subbundle $\nu_g^-(C)$ of $TM|_C$ given by the negative space of the Hessian
of $g$ has constant dimension along $C$, call it ${\rm Ind}C$;
\item the equivariant Euler class $e_T(\nu_g^-(C))$ is not a
zero-divisor in $H^*_T(C)$.
\end{itemize}
This is sufficient to deduce that for any critical set $C$, the
result stated at point a) of the previous theorem remains true if
we replace $f$ by $g$ and $F$ by $C$.

   \section{The residue kernel of the equivariant Kirwan map}

   The goal of this section is to give a proof of Theorem
   \ref{secondmain}.
This theorem is similar to the main result of \cite{Je}: it is a
generalization
of \cite{Je} in that  the proof  does not require the  hypothesis of
isolated fixed points
and the argument is generalized to equivariant Kirwan maps.
Like the argument in \cite{Je}, our argument makes use of residues.

 From the definition of $\ker_{{\rm
   res}}(\kappa_S)$ (see Definition 1.4), it is obvious that the
   latter contains $K_+$. In order to prove that it contains $K_-$ as
   well, we only have to  notice that, by the
   Atiyah-Bott-Berline-Vergne  localization formula (see [At-Bo1] and
   [Be-Ve]), we have that
   $${\rm Res}_X^+\sum_{F\in\cF_+}\int_F\frac{i_F^*(\eta\zeta)}{e_F}=
   -{\rm
   Res}_X^+\sum_{F\in\cF_-}\int_F\frac{i_F^*(\eta\zeta)}{e_F}.$$ The
   difficult part of the proof will be to show that
   $$ \ker_{{\rm res}}(\kappa_S)\subset K_-\oplus K_+.$$

   Let us start with the following two  lemmas.
   \begin{lemma}\label{firstlemma} If $F$ is in $\cF$ and $\alpha\in
 H_T^*(F)$ satisfies
    \begin{equation}\label{hypothesis}\int_F
   i_F^*(\zeta)\alpha=0\end{equation} for all $\zeta \in H^*(M)$, then
   $$\alpha=0.$$
   \end{lemma}
   \begin{proof} We have $$\int_F i_F^*(\zeta)\alpha =\int_M\zeta
 i_{F*}(\alpha),$$
   where $i_{F*}:H_T^*(F)\to H_T^{*+{\rm codim}F}(M)$ is the
   push-forward of the inclusion map $i_F:F\to M$. So  equation
   (\ref{hypothesis}) leads to
   $$i_{F*}(\alpha)=0.$$
   We apply $i_F^*$ on both sides of this equation and  use the fact
   that the map $i_F^*\circ i_{F*}$ is just the multiplication by the
   Euler class $e_F$ in order to deduce that $$\alpha e_F=0.$$ But
   $e_F$ is the same as the product $e(\nu^-F)e(\nu^+F)$, where neither
   factor is a zero divisor  (see Theorem \ref{morse} (b)). So we
   must have $\alpha=0$. \end{proof}

   \begin{lemma}\label{secondlemma}
   Fix $F\in \cF$ and suppose that $\eta\in H_T^*(M)$ satisfies
   \begin{equation} \label{vanish}{\rm
 Res}_X^+\int_F\frac{i_F^*(\eta\beta)}{e(\nu^- F)}=0
   \end{equation} for  all $\beta \in H_T^*(M).$ Then there exists
   $\gamma\in H_T^*(M)$ such that
   $$\eta|_F=\gamma|_F$$
   and $$\gamma|_G=0 \ {\it \ for \ all \ } G\in \cF \ {\rm with } \
   f(G)<f(F).$$
   \end{lemma}
   \begin{proof} Consider a $T$-invariant splitting of the normal bundle
   into a sum of complex line bundles. As explained in Remark 2.1.2,
   this decomposition leads to the series
   $$\frac{i_F^*(\eta)}{e(\nu^- F)}=\sum_{r\le r_0}\gamma_rX^r,$$
   with $\gamma_r\in H^*(F)\otimes \C[\{Y_i\}]$. For $\zeta\in
   H^*(M)$ we have that
   $$\frac{i_F^*(\eta\zeta)}{e(\nu^- F)}=\sum_{r\le
 r_0}i_F^*(\zeta)\gamma_rX^r.$$
   The hypothesis of the lemma implies that
   $$\int_Fi_F^*(\zeta)\gamma_{-1}=0,$$
   for all $\zeta\in H^*(M)$.
   From Lemma \ref{firstlemma} it  follows that $$\gamma_{-1}=0.$$

   Now take  an arbitrary
    class $\zeta$  in $H^*(M)$ and put $\beta=X\zeta$ in (\ref{vanish}). We
    deduce that $\gamma_{-2}=0$. Then we take $\beta=X^2\zeta$ and
    deduce that $\gamma_{-3}=0$ etc.
   Consequently we have that
   $$\eta|_F=i_F^*(\eta)=\alpha e(\nu^-F),$$
   where $$\alpha=\sum_{r=0}^{r_0}X^r\gamma_r\in H_T^*(F).$$ Theorem
   2.2.1 (a) says that the multiplication by $e(\nu^-F)$ is the
   composition of the restriction map $H_T^*(M_+)\to H_T^*(F)$ with
   $r$; so we have $$\eta|_F=r(\gamma')|_F$$ where $\gamma'\in
   H_T^*(M_+)$. Note that by the exact sequence
   $$0\to H_T^{*-{{\rm Ind} F}}(F)\stackrel{r}{\to} H_T^*(M_+)\to
 H_T^*(M_-)\to 0$$
   (see Theorem 2.2.1 (a)) we have that
   $$r(\gamma')|_G=0$$
   for all $G\in \cF$ with $f(G)<f(F)$. Finally we obtain $\gamma\in
   H_T^*(M)$ with the properties required in the lemma by extending
   $r(\gamma')\in H_T^*(M_+)$ to the whole $M$ (see Theorem \ref{morse}
   (c)).  \end{proof}

     We are now ready
   to prove Theorem \ref{secondmain}.

   {\it Proof of Theorem \ref{secondmain}} Denote $$f=\mu_S:M\to
   \s^*=\R$$ the moment map of the $S$ action on $M$. Take $\eta \in
   H_T^*(M)$ satisfying
   \begin{equation}\label{hassumption} {\rm
   Res}_X^+\sum_{F\in{\mathcal
   F}_+}\int_F\frac{i_F^*(\eta\zeta)}{e_F}=0 \end{equation} for all
   $\zeta\in H_T^*(M)$. Consider  the ordering $F_1,F_2, \ldots ,F_N
   $ of the elements of ${\mathcal F}_+$ such that
   $$0<f(F_1)<f(F_2)<\ldots <f(F_N).$$
   We shall inductively construct classes
   $\gamma_1,\gamma_2,\ldots ,\gamma_N\in H_T^*(M)$ which vanish when
   restricted to $\cF_-$ and such that for all $1\leq k \leq N$, the
   form
   $$\eta_k=\eta-\gamma_1- \ldots -\gamma_k $$
   vanishes when restricted to $F_1, F_2, \ldots, F_k $. If we set
   $$\eta_-=\gamma_1+\ldots +\gamma_N,\quad \eta_+=\eta-\eta_-,$$
   then the decomposition $$\eta=\eta_-+\eta_+$$ has the desired
   properties.

   First, in (\ref{hassumption}) we put
   $$\zeta = \beta\alpha^+(F_1),$$ where $\beta\in H_T^*(M)$ is an
   arbitrary cohomology class and the notation
   $\alpha^+(F_1)$ was introduced in Theorem 2.2.1 (d).
     Since
   $$e(\nu_f^-(F_1))e(\nu_f^+(F_1))=e_{F_1},$$
   we deduce that
   $${\rm Res}_{X}^+\int_{F_1}\frac{i_{F_1}^*(\beta\eta)}{e(\nu^-(F_1))}=0
   \ {\rm for \  all \ } \beta \in H_T^*(M),$$ where we are using the
   fact that $\alpha^+(F_1)|_F = 0 $ if $f(F) >
     f(F_1) $. From
   Lemma \ref{secondlemma}  it follows that there exists $\gamma_1\in
   H_T^*(M)$ such that
   $$\eta|_{F_1}=\gamma_1|_{F_1},$$
   and $$\gamma_1|_{G}=0, \ {\rm  for \ all} \ G \in \cF  \ {\rm with
   } \ f(G)<f(F_1).$$

   Next suppose that we have constructed $\gamma_1, \ldots,\gamma_k$
   which vanish when restricted to $\cF_-$ and such that
   $$\eta_k=\eta-\gamma_1-\ldots -\gamma_k$$
   vanishes when restricted to $F_1,\ldots,F_k$. We claim that
   $\eta_k$ satisfies
   $$\text{Res}_X^+\sum_{F\in{\mathcal
   F}_+}\int_F\frac{i_F^*(\eta_k\zeta)}{e_F}=0.$$ This follows from
   equation (\ref{hassumption}) and the Atiyah-Bott-Berline-Vergne
 localization formula for
   $(\gamma_1+\ldots +\gamma_k)\zeta$ (note that the restriction of
   the  form    $(\gamma_1+\ldots +\gamma_k)\zeta$ to ${\mathcal F}_-$ is
 zero). Since $\eta_k|_G=0$
   for all $G\in{\mathcal F}_+$ with $f(G)\leq f(F_k)$, we
   have that \begin{equation}\label{reduced}{\rm
   Res}_X^+\sum_{F\in{\mathcal
   F}_+,f(F) f(F_k)}\int_F\frac{i_F^*(\eta_k\zeta)}{e_F}=0.
   \end{equation} In (\ref{reduced}) we put
   $\zeta=\beta\alpha^+(F_{k+1})$, where $\beta$ is an arbitrary
   element in $H_T^*(M)$. Since
   $$e(\nu_f^-(F_{k+1}))e(\nu_f^+(F_{k+1}))=e_{F_{k+1}},$$
   we deduce that

$$\text{Res}_{X}^+\int_{F_{k+1}}\frac{i_{F_{k+1}}^*(\beta\eta)}{e(\nu^-(F_{k
 +1}))}=0$$
     where we are using the fact that $\alpha^+(F_{k+1})|_F = 0 $ if $f(F) >
     f(F_{k+1}) $. From
   Lemma \ref{secondlemma}  it follows that there exists
   $\gamma_{k+1}\in H_T^*(M)$ such that
   $$\eta_k|_{F_{k+1}}=\gamma_{k+1}|_{F_{k+1}},$$
   and $$\gamma_{k+1}|_G=0$$ for all $G\in \cF$ with $f(G)<
   f(F_{k+1})$. This means that the form
   $$\eta_{k+1}=\eta_k-\gamma_{k+1}=\eta-\gamma_1
   -\ldots -\gamma_{k+1}$$ vanishes  when restricted to $F_{k+1}$ and
   for all $1\leq j\leq k$, we have that
   $$\eta_{k+1}|_{F_j}=\eta_k|_{F_j}=0.$$ \qed

   \section{The kernel of $\kappa:H^*_T(M)\to H^*(M\quott T)$ via
equivariant
 Kirwan maps }

   We will provide a proof of Theorem \ref{firstmain}. A similar result is
stated in [Go] (where it is proved by expressing $\kappa$ as a
composition of  $\kappa_S$). The general strategy
 we will use is the same
   as in the proof of Theorem 3 of [To-We].

   {\it Proof of Theorem \ref{firstmain}} First, it is simple to prove that
 for every generic circle
   $S\subset T$, we have that
   $$\ker\kappa_S\subset \ker \kappa.$$
   We only have to take into account the residue formulas (\ref{jkformula})
 and (\ref{abelianresidue}),
   where in the latter formula res is an iterated residue which starts with
 ${\rm Res}^+_{X}.$

   Let us prove that
   $$\ker\kappa \subset \sum_S\ker\kappa_S.$$
   To this end, we consider the ordering $\mu^{-1}(0)=C_0, C_1, C_2,
   \ldots, C_p$ of the critical sets of $g$
   such that $$i<j \Rightarrow g(C_i)<g(C_j).$$
   Take $\eta$  in $\ker\kappa$, i.e. $\eta|_{C_0}=0$. We will  show by
 induction
   on $0\le k \le p$ that there exists $\eta_k\in \sum_S\ker\kappa_S$ such
 that
   $$\eta_k|_{C_i} = \eta|_{C_i}, \quad {\rm for \ all \ } 0\le i \le k.$$
   We start with $\eta_0=0$, and at the end of the process we will obtain
the
 form $\eta_p\in \sum_S\ker\kappa_S$
 which, by the Kirwan injectivity theorem\footnote{Note that $M^{T}\subset {\rm Crit}(g)$.},
  is the same as $\eta$.

   Suppose that we have  $\eta_k$ and want to construct $\eta_{k+1}$.
   The form $\eta-\eta_k$ vanishes on all critical sets $C$ with
   $$g(C) <g(C{_{k+1}}).$$
  By Remark 2.2.2,
 $(\eta-\eta_k)|_{C_{k+1}}$ is a multiple of
   $e_T(\nu^-_g(C_{k+1}))$. Consider the function
   $$h(x)=\norm\mu(x) +\lambda \mu(C_{k+1})\norm^2, \quad x\in M,$$
   where $\lambda$ is a sufficiently large positive real number (see
   below). One can easily check that $C_{k+1}$ is a critical set
   of $g$ and the negative spaces  of the Hessians of $f$ and $g$
   at any point in $C_{k+1}$ are the same. We use Remark 2.2.2 for
   the function $h$, which is the norm squared of a moment map:
   Because $(\eta -\eta_{k})_{C_{k+1}}$ is a multiple of
   $e_T(\nu_g^-(C_{k+1}))=e_T(\nu_h^-(C_{k+1}))$, we deduce that
   there exists a form
   $\beta_{k}\in H_T^*(M)$ which vanishes on
   $h^{-1}(-\infty, h(C_{k+1})-\epsilon)$ and such that
    $$\beta_{k}|_{C_{k+1}} = (\eta-\eta_{k})|_{C_{k+1}}.$$
   In particular we have that
   \begin{equation} \label{zero}\beta_k|_{C} =0,\quad {\rm for \ all \ critical \ sets  \
 \ } C \ {\rm of \ } \ g {\rm \  with\ } h(C)<h(C_{k+1}).\end{equation}
 So if we set $\eta_{k+1}:=\eta_k+\beta_k$ then obviously
$$(\eta-\eta_{k+1})|_{C_j}=0, \quad {\rm for \ all \ } 0\le j\le k+1.$$

Now we claim that there exists $\lambda\in \R$ such that the form
$\beta_k$ is in $ \sum_S\ker{\kappa_S}$. In fact we only need to
choose $\lambda$ sufficiently large that $h(F)<h(C_{k+1})$ for any
$F\in {\mathcal F}$ with $\langle \mu(F), \mu(C_{k+1})\rangle <0$.
Note that any such $F$ is contained in a critical set $C$, because
$M^T\subset {\rm Crit}(g)$. Consequently $h(C)<h(C_{k+1})$, thus,
by (\ref{zero}), we have that $$\beta_k|_F=0 \quad {\rm for \ all}
\ F\in {\mathcal F}{\rm  \  with \ } \langle \mu(F),
\mu(C_{k+1})\rangle <0.$$ We deduce that there exists a generic
circle $S$ such that $\beta_k\in \ker\kappa_S$. Indeed, the
generic circle $S\subset T$ can be chosen such that
 the two components
   $\mu_S$ and $\langle \mu(\cdot), \mu(C_{k+1})\rangle$  of $\mu$ are
 sufficiently close such that
   $$\mu_S(F)<0 \Leftrightarrow \langle \mu(F),\mu(C_{k+1})\rangle < 0,$$
   where $F\in \cF$. We also use the characterization of
   $\ker\kappa_S$ given by Theorem \ref{secondmain}
   (note that, in the
   notation of Theorem \ref{secondmain}, we have  $\beta_k \in K_-$).
This finishes the proof. \hfill $\square$

   \section{Hamiltonian actions of non-abelian Lie groups}

   The goal of this section is to prove Theorem \ref{nonabelian}. We
   will use the notations from the introduction. The action of $W$ on
   $H^*_T(M)$ plays an important role. We say that an element $\eta$
   of $H^*_T(M)$ is {\it anti-invariant} if
   $$w.\eta=\eps(w)\eta$$
   for all $w\in W$, where $\eps(w)$ denotes the
   signature\footnote{By definition, $\eps(w)$ is $(-1)^{l(w)}$,
   where $l(w)$ is the length of $w$ with respect to the generating
   set of $W$ consisting of the reflections into the walls of a
   fixed  Weyl chamber in $\t$.}  of $w$. Note that
   \begin{equation}\label{sign}\eps(vw)=\eps(v)\eps(w)\end{equation}
   for all $v,w\in W$. One can easily see that the element $\D$ of
   $H^*_T({\rm pt})=S(\t^*)$ which is obtained by multiplying all
   positive roots (see also (\ref{nonabelianresidue}) and
   (\ref{abelianresidue}) from the introduction) is anti-invariant.
   The following result is proved in [Br]:

   \begin{lemma}\label{Brion} {\rm (see [Br])}The set of anti-invariant
 elements
   in $H^*_T(M)$ is $\D \cdot H^*_T(M)^W$.
   \end{lemma}

   Now we are ready to prove Theorem \ref{nonabelian}.

   {\it Proof of Theorem \ref{nonabelian}}. Let us consider the
   pairing $\langle \ , \ \rangle$ on $H^*_T(M)$ given by $$\langle
   \eta,\zeta\rangle =\kappa_T(\eta\zeta)[M\quott T],$$
   $\eta,\zeta\in H^*_T(M)$. The kernel of $\langle \ , \
   \rangle$  is just $\ker(\kappa_T)$. Since the map $$\kappa_T:H^*_T(M)\to
   H^*(M\quott T)$$ is $W$-equivariant, the pairing $\langle \ , \
   \rangle$ is $W$-invariant.

   First we show that (i) is equivalent to
   (iii). Take $\eta \in H^*_T(M)^W$. By equations
(\ref{nonabelianresidue}),
 (\ref{abelianresidue}) and  (\ref{equivalence}),
     the condition
   $$\kappa_K(\eta)=0$$ is equivalent to
   \begin{equation}\label{orthogonal}\langle\D^2 \eta, \zeta\rangle=0, \
{\rm
 for \ all \  \zeta\in} H^*_T(M)^W \end{equation}
    Now we show that (\ref{orthogonal}) is equivalent to $\D^2\eta \in \ker
 \langle \ , \ \rangle$.
   Indeed, if $\zeta$ is an  arbitrary element in $H^*_T(M)$, we can
consider
   $$\zeta'=\sum_{w\in W}w\zeta$$
   which is $W$-invariant, hence
   \begin{equation}\label{finish}\langle \D^2\eta,\zeta'\rangle=|W|\langle
 \D^2\eta,\zeta\rangle\end{equation}
   where we have used the $W$-invariance of both the pairing $\langle \ , \
 \rangle$ and the
   class $\D^2\eta$. Equation (\ref{finish}) finishes the proof of the
 equivalence between
   (i) and (iii).

   Now we prove that (ii) is equivalent to (iii). In fact only the
 implication
   $\rm{(iii)\Rightarrow (ii)}$ is non-trivial. So let us consider $\eta\in
 H^*_T(M)^W$
   with the property that $\kappa_T(\D^2\eta)=0$, which is equivalent to
   $$\langle \D^2\eta,\zeta\rangle=0$$
   for every $\zeta\in H^*_T(M)$.
   \newline {\it Claim 1.} $\langle \D\eta,\zeta\rangle=0$, for all
$\zeta\in
 H^*_T(M)$ which
   is anti-invariant.
   \newline Indeed, if $\zeta$ is anti-invariant, by Lemma \ref{Brion}, it
is
 of the form
   $\zeta=\D\xi$, with $\xi\in H^*_T(M)$, thus we have that
   $$\langle \D\eta,\zeta\rangle=\langle \D\eta,\D\xi\rangle=\langle
 \D^2\eta,\xi\rangle=0.$$
   \newline{\it Claim 2.} $\langle \D\eta,\zeta\rangle=0$, for all
$\zeta\in
 H^*_T(M)$.
   \newline This is because we can use Claim 1 for
   $$\zeta'=\sum_{w\in W}\eps(w)w.\zeta$$
   which is an anti-invariant equivariant cohomology class. More precisely,
 we  notice that
   \begin{equation}\label{prime}\langle \D\eta,\zeta'\rangle=|W|\langle
 \D\eta,\zeta\rangle,\end{equation}
   where we have used the $W$-invariance of $\langle \ , \ \rangle$ and the
   anti-invariance of $\D\eta$. The left hand side of (\ref{prime}) is zero
 by Claim 1,
   hence so must be the right hand side.
   $\hfill \square$

   {\bf Remark.} Explicit descriptions of $\ker(\kappa_K)$ in terms of
 $\ker(\kappa_T)$
   can be easily deduced from Theorem \ref{nonabelian}, as follows:
   \begin{align}\label{firstchar}\ker(\kappa_K)&=\{\frac{1}{\D}\sum_{w\in
 W}\eps(w)w.\eta ~|~ \eta \in \ker (\kappa_T)\}
   \end{align}
   Alternatively, we obtain all elements of the kernel of $\kappa_K$
   if we consider all cohomology classes of the type $\sum_{w\in W}w.\eta$
   with $\eta\in \ker(\kappa_T)$ which are multiples of $\D^2$, and divide
 those by $\D^2$
   (of course the description (\ref{firstchar})
   is simpler, since by Lemma \ref{Brion}, for all $\eta\in H^*_T(M)$, the
 cohomology class
    $\sum_{w\in W}\eps(w)w.\eta$ is a multiple of $\D$).
   The presentation (\ref{firstchar}) is also given  in [To-We, Proposition
 6.1] (see also
   [Br, Corollaire 1]).


\appendix 
\section{(By Jonathan Woolf)}

This appendix explains how to prove Theorem \ref{nonabelian} under weaker assumptions. In particular the assumption that $0$ is a regular value for the  $T$  moment map can be relaxed in a wide class of cases. When $0$ is not a regular value for the $T$ moment map $M\quott T$ is, in general, singular. It has a natural stratification by symplectic orbifolds \cite{sl} and, in particular, is a pseudomanifold with even dimensional strata. Thus we can define its intersection cohomology $\ih^*(M \quott T)$ (with coefficients in $\mathbb{C}$) and there is a natural non-degenerate intersection pairing
$$
\langle\  , \ \rangle_{M\quott T} : \ih^*(M \quott T) \times \ih^*(M \quott T) \to \mathbb{C}
$$
generalising the intersection pairing on the cohomology of a manifold. We can extend the definition of the Kirwan map to obtain a map
$$
\kappa_T : H^*_T(M) \to \ih^*(M\quott T).
$$
This construction was first carried out for the closely related case of geometric invariant theory quotients in algebraic geometry in \cite{k3}.  A version for symplectic quotients can be found in \cite{kw1}. It should be noted that when $0$ is not a regular value of $\mu_T$ there are choices involved in this construction, and so $\kappa_T$  is not canonical. We briefly sketch the construction.

Meinrenken and Sjamaar describe a partial desingularisation (orbifold singularities may remain) procedure for singular symplectic quotients in \cite{mes}. (This is modelled on the algebro-geometric version in \cite{k2}.) Singularities in the quotient arise from Lie subgroups of $T$ with fixed point subsets contained within $\mu_T^{-1}(0)$. There is a finite set $\mathcal{S}_M$ of such subgroups with dimension $\geq 1$, and $0$ is a regular value of $\mu_T$ if, and only if, this set is empty. If $S \in \mathcal{S}_M$ is of maximal dimension then the components of its fixed point subset contained within $\mu_T^{-1}(0)$ form a closed symplectic submanifold of $M$. Performing a $T$-equivariant symplectic blowup along this submanifold produces a symplectic manifold $M_1$ with a Hamiltonian $T$ action such that $\mathcal{S}_{M_1} = \mathcal{S}_{ M} - \{S\}$. Proceeding inductively 
we obtain a symplectic manifold $\widetilde M = M_r$ with a Hamiltonian $T$ action and moment map $\tilde\mu_T$ such that $\mathcal{S}_{\widetilde M} = \emptyset$, or equivalently, $0$ is a regular value of $\tilde \mu_T$. The blowdown $M_i \to M_{i-1}$ induces a continuous surjection $\pi_i : M_i\quott T \to M_{i-1}\quott T$. Composing these we obtain a continuous surjection from the symplectic orbifold $\widetilde M \quott T$ onto $M\quott T$ i.e. a partial desingularisation.

For each $i>0$ we can choose, non-canonically,  a surjection $$\ih^*(M_i\quott T) \to \ih^*(M_{i-1}\quott T),$$ see \cite[Theorem 5]{kw1}. We  define a Kirwan map $\kappa_T$ by composing the equivariant pullbacks with  the Kirwan map $\tilde \kappa_T$ and these surjections:
$$
\label{kirwan map}
\xymatrix{
H^*_T(M) \ar@{-->}[d]_{\kappa_T}\ar[r] & H^*_T(M_1) \ar[r] & \ldots  \ar[r] & H^*_T(\widetilde M) \ar[d]^{\tilde \kappa_T} \\
\ih^*(M\quott T) & \ih^*(M_1\quott T)  \ar[l]& \ldots \ar[l]&  \ih^*(\widetilde M\quott T)\cong H^*(\widetilde M\quott T)   \ar[l].
}
$$
What makes the Kirwan map for nonsingular quotients (when $0$ is a regular value of $\mu_T$) useful is its surjectivity. Unfortunately, it is not known whether the map defined above is always surjective when $0$ is not a regular value of $\mu_T$. However, we do have
\begin{theorem}
\label{surjectivity theorem}
Suppose $0$ is not a regular value of $\mu_T$. Then, if the action of $T$ on $M$ is almost-balanced in the sense of \cite[\S 5]{kiem} \emph{or} if $M$ is a complex projective variety with symplectic structure given by the Fubini-Study form and the action of $T$ is the restriction of an algebraic action of $(\mathbb{C}^*)^n$, the Kirwan map $\kappa_T$ is surjective. 
\end{theorem}
\begin{proof}
See \cite{kw2}, \cite{k3} and \cite{w}.\end{proof}

We can use the Kirwan map to define a pairing on $H^*_T(M)$ by
$$
\langle\eta , \zeta \rangle = \langle\kappa_T(\eta) , \kappa_T(\zeta) \rangle_{M\quott T}
$$
where the RHS is the intersection pairing on $\ih^*(M\quott T)$. This extends the previous definition since $\ih^*(M\quott T)\cong H^*(M\quott T)$ and $\langle\kappa_T(\eta) , \kappa_T(\zeta) \rangle_{M\quott T} =  \kappa_T(\eta\zeta) [M\quott T]$ when $0$ is a regular value. If the Kirwan map is surjective we have 
\begin{equation}
\label{ker kappa}
\kappa_T(\eta)=0 \iff \langle \eta, \zeta\rangle = 0 \qquad \forall \zeta \in H^*_T(M).
\end{equation}

Suppose the Hamiltonian $T$ action on $M$ arises as the restriction of a Hamiltonian $K$ action. Then the normaliser $N_K(T)$ acts and preserves $\mu_T^{-1}(0)$, thereby inducing an action of the Weyl group $W = N_K(T)/T$ on $M \quott T$. 
\begin{lemma}
The Kirwan map $\kappa_T : H^*_T(M) \to \ih^*(M\quott T)$ can be chosen to be $W$-equivariant.
\end{lemma}
\noindent  {\em Proof:}
The symplectic blowups in the partial desingularisation procedure can be done $N_K(T)$-equivariantly (so that $N_K(T)$ acts on each $M_i$ and $W$ on each $M_i \quott T$). The Kirwan map $\tilde \kappa_T$ is $W$-equivariant. It remains to see that the surjections $\ih^*(M_i\quott T) \to \ih^*(M_{i-1}\quott T)$ can be chosen $W$-equivariantly. These surjections are constructed in \cite[\S 2]{kw1} by decomposing appropriate complexes of sheaves in the constructible derived category of $M_{i-1}\quott T$. In the presence of a finite group action we can carry out a formally identical argument in the equivariant constructible derived category to obtain the required equivariant surjections.
\hfill $\square$
\begin{corollary}
The pairing $\langle\ ,\ \rangle$ on $H^*_T(M)$ is $W$-invariant.
\end{corollary}
\begin{proof}
This follows immediately from the $W$-equivariance of $\kappa_T$ and the $W$-invariance of the pairing $\langle\ ,\ \rangle_{M \quott T}$ on $\ih^*(M\quott T)$. 
\end{proof}

Now suppose further that $0$ is a regular value for the $K$ moment map $\mu_K$ on $M$. Then $\mu_K^{-1}(0)$ is a submanifold of $M$ and the inclusion $\mu_K^{-1}(0) \hookrightarrow  \mu_T^{-1}(0)$ has trivial normal bundle of dimension $d = \dim \mathfrak{k}^*/\mathfrak{t}^*$. However it is not \emph{equivariantly} trivial --- the equivariant Thom class is the product $\mathcal{D}$ of the positive roots. Provided that no finite subgroups of $T$ fix points of  $\mu_K^{-1}(0)$ then the inclusion $ \mu_K^{-1}(0) / T \hookrightarrow \mu_T^{-1}(0) / T = M\quott T$ is also normally nonsingular. 

Recall from \cite[\S 5.4]{gm2} that, if $\imath : X \hookrightarrow Y$ is a normally nonsingular inclusion of pseudomanifolds with even dimensional stratifications, then there are maps
$$
\imath^* : \ih^*(Y) \to \ih^*(X) \ \textrm{and} \ \imath_*:  \ih^*(X) \to \ih^{*+d}(Y).
$$
These are adjoint to one another with respect to the intersection pairings i.e.
$$
\langle \alpha , \imath_* \beta \rangle_Y = \langle \imath^* \alpha , \beta \rangle_X.
$$
\begin{lemma}
\label{kappa lemma}
Suppose that no finite subgroups of $T$ fix points of  $\mu_K^{-1}(0)$. Let  $\imath : \mu_K^{-1}(0) / T \hookrightarrow M\quott T$ be the normally nonsingular inclusion. Then $\kappa_T(\mathcal{D}\eta) = \imath_*\imath^*\kappa_T(\eta)$.
\end{lemma}
\begin{proof}
For the case when $0$ is regular for $\mu_T$ this reduces to the statement that the 
identification $H^*_T(\mu_T^{-1}(0)) \cong H^*(\mu_T^{-1}(0)/T)$ takes the equivariant Thom class 
$\mathcal{D}$ of the submanifold $\mu_K^{-1}(0)$ to the Thom class of the quotient $\mu_K^{-1}(0)/T$ 
in $\mu_T^{-1}(0)/T$. Moreover this statement is local to $\mu_K^{-1}(0)$. Hence the identity holds more 
generally because, when we partially desingularise, we only blow up along centres which do not meet 
$\mu_K^{-1}(0)$ i.e. we do not alter anything in a neighbourhood of $\mu_K^{-1}(0)$.\end{proof}
 
{\bf Remark.} The restriction  that no finite subgroups of $T$ fix points of  $\mu_K^{-1}(0)$ is not necessary for the proof of Lemma \ref{kappa lemma} (and hence for what follows). If we relax it then the inclusion $\mu_K^{-1}(0) / T \hookrightarrow M\quott T$ will not necessarily be normally nonsingular, the normal fibres may be orbifolds. However, the extensions to Goresky and MacPherson's results on normally nonsingular inclusions can be obtained by working in an appropriate equivariant derived category cf. \cite[Remark 5 and \S 3]{kw1}.

\begin{corollary}
\label{final cor}
For $\eta,\zeta \in H^*_T(M)$ we have
$
\langle \mathcal{D}\eta , \mathcal{D}\zeta\rangle  = \langle \mathcal{D}^2\eta , \zeta\rangle
$.
If both $\eta$ and $\zeta$ are invariant under the Weyl group action i.e.  $\eta,\zeta \in H^*_T(M)^W \cong H^*_K(M)$ then
$$
\langle \mathcal{D}\eta , \mathcal{D}\zeta\rangle  =  \frac{a}{b}|W| \cdot \kappa_K(\eta\zeta)[M\quott K]
$$
where $a$ is the order of the largest subgroup of $K$ which fixes $\mu_K^{-1}(0)$ point-wise and $b$ the order of the largest subgroup of $T$ which fixes $\mu_T^{-1}(0)$ point-wise. 
\end{corollary}
\noindent {\em Proof:}
Using Lemma \ref{kappa lemma} and the fact that $\imath_*$ is adjoint to $\imath^*$ we obtain
$$
\langle \mathcal{D}\eta , \mathcal{D}\zeta\rangle =  \langle \imath^*\kappa_T(\mathcal{D}\eta) , \imath^*\kappa_T(\zeta) \rangle_{\mu_K^{-1}(0) / T}
$$
A second use of Lemma \ref{kappa lemma} and adjointness yields 
$$
\langle  \imath^*\kappa_T(\mathcal{D}\eta) ,\imath^*\kappa_T(\zeta)\rangle_{\mu_K^{-1}(0) / T} = 
\langle\mathcal{D}^2 \eta , \zeta\rangle. 
$$
If $\eta,\zeta \in H^*_T(M)^W$ then, since $0$ is a regular value of $\mu_K$, Theorem B$''$ of \cite{martin} tells us that
$$
 \langle  \imath^*\kappa_T(\mathcal{D}\eta) ,\imath^*\kappa_T(\zeta)\rangle_{\mu_K^{-1}(0) / T} = \frac{a}{b} |W| \cdot \kappa_K(\eta\zeta)[M\quott K].
$$
\hfill$\square$

\begin{theorem}
Suppose $0$ is a regular value of $\mu_K$. Choose a $W$-equivariant Kirwan map $\kappa_T : H^*_T(M) \to \ih^*(M\quott T)$ and suppose that it is surjective. For $\eta \in H^*_T(M)^W \cong H^*_K(M)$ the following assertions are equivalent:
   \begin{itemize}
   \item[(i)] $\kappa_K(\eta)=0$,
   \item[(ii)] $\kappa_T(\mathcal{D}\eta)=0$,
   \item[(iii)] $\kappa_T(\mathcal{D}^2\eta)=0$.
   \end{itemize}
\end{theorem}
\begin{proof}
 Note that $\kappa_K(\eta)=0$ if, and only if, $\kappa_K(\eta\zeta)[M\quott K]=0$ for all $\zeta \in H^*_T(M)^W$. By Corollary \ref{final cor} this occurs if, and only if, $\langle \mathcal{D}\eta , \mathcal{D}\zeta\rangle = 0$ for all $\zeta \in H^*_T(M)^W$. Recall from Lemma \ref{Brion} (which does not require the hypothesis that $0$ is a regular value of $\mu_T$) that
$$
\{ \mathcal{D}\zeta : \zeta \in H^*_T(M)^W \} = \{ \textrm{Anti-invariant}\ \xi \in H^*_T(M)\}.
$$
So $\langle \mathcal{D}\eta , \xi\rangle = 0$ for all anti-invariant $\xi \in H^*_T(M)$ and hence, by the 
argument in the proof of Theorem \ref{nonabelian}, for all $\xi \in H^*_T(M)$. By (\ref{ker kappa}) we see that 
$\kappa_T(\mathcal{D}\eta)=0$ i.e. (i) implies (ii). The converse is easy. 
 
On the other hand, using Corollary \ref{final cor} we see that $\kappa_K(\eta)=0$ if, and only if,
 $\langle \mathcal{D}^2\eta , \zeta\rangle = 0$ for all $\zeta \in H^*_T(M)^W$. Since $\mathcal{D}^2\eta$ is 
$W$-invariant this is equivalent to $\langle \mathcal{D}^2\eta , \xi\rangle = 0$ for all $\xi \in H^*_T(M)$. 
By (\ref{ker kappa}) this  is equivalent to $\kappa_T(\mathcal{D}^2\eta)=0$. So (i) and (iii) are equivalent. 
\end{proof}

In particular this theorem applies under the conditions in the statement of Theorem \ref{surjectivity theorem}. It can be read either as allowing us to determine the kernel of $\kappa_K$ in terms of the kernel of {\em any} surjective $W$-equivariant Kirwan map $\kappa_T$ for the $T$ action, or, alternatively, as telling us that the subspace $\{ \mathcal{D}\eta : \eta \in \ker(\kappa_K)\}$ is contained in the kernel of {\em every}  surjective $W$-equivariant Kirwan map $\kappa_T$.


   \bibliographystyle{abbrv}

\begin{thebibliography}{Kie-Wo2}

   \bibitem[At-Bo1]{At-Bo1}  M. F. Atiyah and R. Bott, {The moment map and
   equivariant cohomology}, {\em Topology}
    {\bf 23} (1984),  1--28.

   \bibitem[At-Bo2]{At-Bo2}  M. F. Atiyah and R. Bott,  {The Yang-Mills
 equations
   over Riemann surfaces}, {\em Phil. Trans. Roy. Soc. London}   {\bf
   A308} (1982), 523-615.

   \bibitem[Be-Ve]{Be-Ve}  N. Berline  and M. Vergne,  {Z\'eros d'un champ
   de vecteurs et classes caract\'eristiques \'equivariantes},  {\em Duke
   Math. J.}
    {\bf  50}  (1983), 539--549.

   \bibitem[Bo-Tu]{Bo-Tu} R. Bott and L. Tu, {\em Differential Forms in
 Algebraic Topology},
    Springer-Verlag, 1982.

   \bibitem[Br]{Br} M. Brion, {Cohomologie \'equivariante des points
 semi-stables},
   {\em J. reine angew. Math.}, {\bf 421} (1991), 125-140.




   \bibitem[Go]{Go} R. F. Goldin, {An effective algorithm
   for the cohomology ring of symplectic reductions},
   {\em Geom.\ and Func.\
   Anal.}, {\bf 12} (2002), 567--583.

\bibitem[Go-Ma]{gm2} M. Goresky and R. MacPherson, {Intersection
Homology Theory {II}}, {\em Invent. Math.} {\bf 71} (1983), 77--129. 

   \bibitem[Gu-St]{Gu-St} V. Guillemin and S. Sternberg, {\em Supersymmetry
 and Equivariant de Rham
   Theory}, Springer-Verlag, 1999.

 \bibitem[Gu-K]{Gu-K} V. Guillemin and J. Kalkman,
 The Jeffrey-Kirwan theorem and residue operations in
 equivariant cohomology, {\em Jour. reine angew. Math.}
{\bf   458} (1995), 37--52.


   \bibitem[Je]{Je}  L. C. Jeffrey, { The residue formula and the
   Tolman-Weitsman theorem}, {\em Jour. reine  angew. Math.} 
{\bf 562} (2003), 51--58.


   \bibitem[Je-Ki1]{Je-Ki1} L. C. Jeffrey and F. C. Kirwan, {
   Localization for nonabelian group actions}, {\em Topology}
      {\bf 34} (1995), 291--327.



   \bibitem[Je-Ki2]{Je-Ki2} L. C. Jeffrey and F. C. Kirwan, {
   Localization and the quantization conjecture}, {\em Topology}
    {\bf 36} (1995),               647--693.



   \bibitem[Je-Ma]{Je-Ma}
   L. C. Jeffrey and A.-L. Mare, {The kernel of the equivariant
   Kirwan map and the residue formula}, {\it Quart. J. Math. Oxford}, {\bf 54} (2004), 435--444.

\bibitem[Kie]{kiem} Y-H. Kiem, {Intersection cohomology of quotients of nonsingular
varieties}, {\em Invent. Math.} {\bf 155} No. 1 (2004), 163--202.

\bibitem[Kie-Wo1]{kw1} Y-H. Kiem and J. Woolf, {The {K}irwan map for singular symplectic quotients}, {\em Preprint} 2004.

\bibitem[Kie-Wo2]{kw2} Y-H. Kiem and J. Woolf, 
{The cosupport axiom, equivariant cohomology and the intersection cohomology of certain singular symplectic quotients}, 
preprint {\tt math.AG/0101255}.



   \bibitem[Ki]{Ki}  F. C. Kirwan, {\em Cohomology of Quotients in
   Complex and Algebraic Geometry},  Mathematical Notes {\bf 31},
   Princeton University Press, Princeton N. J., 1984.


\bibitem[Ki2]{k2} F. Kirwan,
{Partial Desingularisations of
Quotients of Nonsingular Varieties and their {B}etti Numbers},
{\em Annals of Mathematics} {\bf122} (1985), 41--85.

\bibitem[Ki3]{k3} F. Kirwan,
{Rational Intersection Homology of
Quotient Varieties}, {\em Invent. Math.}  {\bf 86} (1986), 471--505.

\bibitem[Ma]{martin} S. Martin, {Symplectic quotients by a nonabelian group and by its maximal torus},
{\em Annals of Math.} (to appear), preprint {\tt math.SG/0001001}.

\bibitem[Me-Sj]{mes} E. Meinrenken and R. Sjamaar, {Singular reduction and quantization},
 {\em Topology}  {\bf 38} No. 4 (1999), {699--762}.

\bibitem[Sj-Le]{sl}
R. Sjamaar and E. Lerman, {Stratified Symplectic Spaces and Reduction},
{\em Annals of Math.}  {\bf 134} (1991), 375--422.




   \bibitem[To-We]{To-We}  S. Tolman and J. Weitsman,  The cohomology rings of symplectic quotients,
 {\it Comm. Anal. Geom.}  {\bf 11}  No. 4 (2003), 751–-773.


\bibitem[Wo]{w}
J. Woolf, {The decomposition theorem and the intersection cohomology of quotients in algebraic geometry},
{\em J. Pure and Applied Algebra}  {\bf 182} (2003), 317--328.





   \end{thebibliography}

  \end{document}